
\input amstex
\documentstyle{amsppt}


\expandafter\ifx\csname amsppt.sty\endcsname\relax\input amsppt.sty \fi

\chardef\tempatcatcode\the\catcode`\@
\catcode`\@=11

\W@{This is LABEL.DEF by A.Degtyarev <March 9, 1997>}
\ifx\labelloaded@\undefined\else
	\catcode`\@\tempatcatcode\let\tempatcatcode\undefined
  \message{[already loaded]}\endinput\fi
\let\labelloaded@\tempatcatcode\let\tempatcatcode\undefined
\def\labelmesg@ {LABEL.DEF: }

\def\stylefile@#1{\expandafter\ifx\csname#1\endcsname\relax
	\else\message{[already loaded]}\endinput\fi
	\expandafter\edef\csname#1\endcsname{\catcode`\noexpand\@\tempcat
		\toks@{}\toks@@{}\expandafter\let\csname#1\endcsname\noexpand\empty}%
	\let\tempcat\undefined\uppercase{\def\styname{#1}}}
\def\loadstyle#1{\edef\next{#1}%
	\DN@##1.##2\@nil{\if\notempty{##2}\else\def\next{##1.sty}\fi}%
	\expandafter\next@\next.\@nil
	\expandafter\ifx\csname\next\endcsname\relax\input\next\fi}

\def\make@letter{\edef\t@mpcat{\catcode`\@\the\catcode`\@}\catcode`\@=11 }

\def\arabicnum#1{\number#1}

\def\Romannum#1{\expandafter\uppercase\expandafter{\romannumeral #1}}
\def\alphnum#1{\ifcase#1\or a\or b\or c\or d\else\@ialph{#1}\fi}
\def\@ialph#1{\ifcase#1\or \or \or \or \or e\or f\or g\or h\or i\or j\or
	k\or l\or m\or n\or o\or p\or q\or r\or s\or t\or u\or v\or w\or x\or y\or
	z\else\fi}
\def\Alphnum#1{\ifcase#1\or A\or B\or C\or D\else\@Ialph{#1}\fi}
\def\@Ialph#1{\ifcase#1\or \or \or \or \or E\or F\or G\or H\or I\or J\or
	K\or L\or M\or N\or O\or P\or Q\or R\or S\or T\or U\or V\or W\or X\or Y\or
	Z\else\fi}

\def\@car#1#2\@nil{#1}
\def\@cdr#1#2\@nil{#2}
\newskip\@savsk
\let\@ignorespaces\ignorespaces
\def\ignorespaces{\ifhmode
  \ifdim\lastskip>\z@\else\nobreak\hskip1sp minus1sp\fi\fi\@ignorespaces}
\def\@bsphack{\relax\ifmmode\else\@savsk\lastskip
  \ifhmode\edef\@sf{\spacefactor\the\spacefactor}\fi\fi}
\def\@esphack{\relax
  \ifx\penalty@\penalty\else\penalty\@M\fi   
  \ifmmode\else\ifhmode\@sf{}\ifdim\@savsk>\z@\@ignorespaces\fi\fi\fi}
\newread\@inputcheck
\def\@input#1{\openin\@inputcheck #1 \ifeof\@inputcheck \W@
  {No file `#1'.}\else\closein\@inputcheck \relax\input #1 \fi}

\def\eat@bs{\expandafter\eat@\string}
\def\eat@ii#1#2{}
\def\eat@iii#1#2#3{}
\def\eat@iv#1#2#3#4{}
\def\@xname#1{\expandafter\noexpand\csname\eat@bs#1\endcsname}
\def\@name#1{\csname\eat@bs#1\endcsname}
\def\@edefname#1{\expandafter\edef\csname\eat@bs#1\endcsname}
\def\@xdefname{\global\@edefname}
\long\def\@ifundefined#1#2#3{\expandafter\ifx\csname\eat@bs#1\endcsname\relax
  #2\else#3\fi}
\def\@@addto#1#2{{\toks@\expandafter{#1#2}\xdef#1{\the\toks@}}}
\def\@@addparm#1#2{{\toks@\expandafter{#1{##1}#2}%
	\edef#1{\gdef\noexpand#1####1{\the\toks@}}#1}}
\def\ST@P{@step}
\def\ST@LE{style}
\def\N@M{no}
\outer\def\newcounter{\checkbrack@{\expandafter\newcounter@\@txtopt@{{}}}}
{\let\newcount\relax
\gdef\newcounter@#1#2#3{{
	\toks@@\expandafter{\csname\eat@bs#2\N@M\endcsname}%
	\DN@{\alloc@0\count\countdef\insc@unt}%
	\ifx\@txtopt@\identity@\expandafter\next@\the\toks@@
		\else\if\notempty{#1}
			\global\expandafter\let\csname\eat@bs#2\N@M\endcsname#1\fi\fi
  \@ifundefined{\the\eat@bs#3}{\toks@{}}{%
		\toks@\expandafter{\csname the\eat@bs#3\endcsname.}}%
	\@xdefname{\the\eat@bs#2}{\the\toks@\noexpand\arabicnum\the\toks@@}%
  \@xdefname{#2\ST@P}{}%
  \@ifundefined{#3\ST@P}{}
  	{\edef\next@{\noexpand\@@addto\@xname{#3\ST@P}{%
			 \global\@xname{#2\N@M}\z@\@xname{#2\ST@P}}}\next@}%
	\expandafter\@@addto\expandafter\closeaux@\expandafter
		{\expandafter\\\the\toks@@}}}}
\outer\def\copycounter#1#2{%
	\@xdefname{#1\N@M}{\@xname{#2\N@M}}%
	\@xdefname{#1\ST@P}{\@xname{#2\ST@P}}%
	\@xdefname{\the\eat@bs#1}{\@xname{\the\eat@bs#2}}}
\outer\def\everyreset{\checkstar@\everyreset@}
\def\everyreset@#1{\ifx\@numopt@\identity@\let\next@\@@addto
  \else\let\next@\gdef\fi\expandafter\next@\csname\eat@bs#1\ST@P\endcsname}
\def\counterstyle#1{\expandafter\gdef\csname the\eat@bs#1\endcsname}
\def\advancecounter#1#2{\@name{#1\ST@P}\global\advance\@name{#1\N@M}#2}
\def\setcounter#1#2{\@name{#1\ST@P}\global\@name{#1\N@M}#2}
\def\counter#1{\refstepcounter#1\printcounter#1}
\def\printcounter#1{\@name{\the\eat@bs#1}}
\def\refcounter#1{\xdef\@lastmark{\printcounter#1}}
\def\stepcounter#1{\advancecounter#1\@ne}
\def\refstepcounter#1{\advancecounter#1\@ne\refcounter#1}
\def\savecounter#1{\@edefname{#1@sav}{%
	\global\@name{#1\N@M}\the\@name{#1\N@M}}}
\def\restorecounter#1{\@ifundefined{#1@sav}{}{\@name{#1@sav}}}

\def\warning#1#2{\W@{Warning: #1 on input line #2}}
\def\warning@#1{\warning{#1}{\the\inputlineno}}
\def\bftext{\ifmmode\fam\bffam\else\bf\fi}
\let\@lastmark\empty
\let\@lastlabel\empty
\def\lastmark{\@lastmark}
\let\lastlabel\empty
\def\newlabel#1#2#3{{\@ifundefined{\r@-#1}{}{\warning
  {label `#1' multiply defined}{#3}}%
  {\let\protect\noexpand\@xdefname{\r@-#1}{#2}}}}
\def\@@@xref#1{%
	\@ifundefined{\r@-#1}{{\bftext??}\warning@{label `#1' undefined}}%
  {\expandafter\expandafter\expandafter\@car\csname r@-#1\endcsname\@nil}}
\def\@xref#1{\rom{\@@@xref{#1}}}
\let\xref\@xref
\def\pageref#1{%
	\@ifundefined{\r@-#1}{{\bftext??}\warning@{label `#1' undefined}}%
  {\expandafter\expandafter\expandafter\@cdr\csname r@-#1\endcsname\@nil}}
\def\thepage{\ifnum\pageno<\z@\romannumeral-\pageno\else\number\pageno\fi}
\def\label@#1#2{\@bsphack\edef\@lastlabel{#1}{\let\thepage\relax
  \def\protect{\noexpand\noexpand\noexpand}%
  \edef\@tempa{\write\@auxout{\string
    \newlabel{#2}{{\@lastmark}{\thepage}}{\the\inputlineno}}}%
  \expandafter}\@tempa\@esphack}
\def\label#1{\label@{#1}{#1}}
\def\@gobble{\relaxnext@
 	\DN@{\ifx[\next\DN@[####1]{}\else
 		\ifx"\next\DN@"####1"{}\else\DN@{}\fi\fi\next@}%
 	\FN@\next@}
\def\@gobblefn#1{\ifx#1[\expandafter\@gobblebr\else
  \ifx#1"\expandafter\expandafter\expandafter\@gobblequ\fi\fi}
\def\@gobblebr#1]#2{}
\def\@gobblequ#1"#2{}
{\catcode`\~\active\lccode`\~=`\@
\lowercase{\global\let\save@at=~ \gdef\protect@at{\def~{\protect\save@at}}}}
\def\Protect@@#1{\def#1{\protect#1}}
\def\disable@special{\let\@writeaux\eat@ii\let\label\eat@
	\def\footnotemark{\protect\@gobble}%
  \let\footnotetext\@gobblefn\let\footnote\@gobblefn
  \Protect@@\\\let\ifvmode\iffalse
	\let\refcounter\eat@\let\advancecounter\eat@ii\let\setcounter\eat@ii}
\def\@writeaux#1#2{\@bsphack{\disable@special\protect@at
	\def\chapter{\protect\chapter@toc}\let\thepage\relax
	\def\protect{\noexpand\noexpand\noexpand}%
	\Protect@@~\Protect@@\@@@xref\Protect@@\pageref\Protect@@\nofrills
  \edef\@tempa{\@ifundefined#1{}{\write#1{#2}}}\expandafter}\@tempa\@esphack}
\def\writeauxline#1#2#3{\@writeaux\@auxout
  {\string\@auxline{#1}{#2}{#3}{\thepage}}}
{\let\newwrite\relax
\gdef\@openin#1{\make@letter\@input{\jobname.#1}\t@mpcat}
\gdef\@openout#1{\global\expandafter\newwrite\csname tf@-#1\endcsname
   \immediate\openout\csname tf@-#1\endcsname \jobname.#1\relax}}
\def\auxlinedef#1{\expandafter\def\csname do@-#1\endcsname}
\def\@auxline#1{\@ifundefined{\do@-#1}{\expandafter\eat@iii}%
	{\expandafter\expandafter\csname do@-#1\endcsname}}
\def\begin@write#1{\bgroup\def\do##1{\catcode`##1=12 }\dospecials\do\~\do\@
	\catcode`\{=\@ne\catcode`\}=\tw@\immediate\write\csname#1\endcsname}
\def\end@writetoc#1#2#3{{\string\tocline{#1}{#2\string\page{#3}}}\egroup}
\def\do@tocline#1{%
	\@ifundefined{\tf@-#1}{\expandafter\eat@iii}
		{\begin@write{tf@-#1}\expandafter\end@writetoc}
}
\auxlinedef{toc}{\do@tocline{toc}}

\let\protect\empty
\def\Protect#1{\@ifundefined{#1@P@}{\PROTECT#1}{}}
\def\PROTECT#1{%
	\expandafter\let\csname\eat@bs#1@P@\endcsname#1%
	\edef#1{\noexpand\protect\@xname{#1@P@}}}

\Protect\operatorname
\Protect\operatornamewithlimits
\Protect\qopname@
\Protect\qopnamewl@
\Protect\text
\Protect\topsmash
\Protect\botsmash
\Protect\smash
\Protect\widetilde
\Protect\widehat
\Protect\thetag
\Protect\therosteritem
\Protect\Cal
\Protect\Bbb
\Protect\bold
\Protect\slanted
\Protect\roman
\Protect\italic
\Protect\boldkey
\Protect\boldsymbol
\Protect\frak
\Protect\goth
\Protect\dots
\Protect\cong
\Protect\lbrace \let\{\lbrace
\Protect\rbrace \let\}\rbrace
\let\root@P@@\root \def\root@P@#1{\root@P@@#1\of}
\def\root#1\of{\protect\root@P@{#1}}

\let\@frills@\identity@
\let\@txtopt@\identyty@
\let\@numopt@\identity@
\def\frills{\ignorespaces\@txtopt@}
\def\numberline{\@numopt@}
\newif\if@write\@writetrue
\def\checkstar@#1{\DN@{\@writetrue
  \ifx\next*\DN@####1{\let\@numopt@\eat@\checkstar@@{#1}}%
	  \else\DN@{\let\@numopt@\identity@#1}\fi\next@}\FN@\next@}
\def\checkstar@@#1{\DN@{%
  \ifx\next*\DN@####1{\@writefalse#1}%
	  \else\DN@{\@writetrue#1}\fi\next@}\FN@\next@}
\def\checkfrills@#1{\DN@{%
  \ifx\next\nofrills\DN@####1{#1}\def\@frills@####1{####1\nofrills}%
	  \else\DN@{#1}\let\@frills@\identity@\fi\next@}\FN@\next@}
\def\checkbrack@#1{\DN@{%
	\ifx\next[\DN@[####1]{\def\@txtopt@########1{####1}%
		#1}%
	\else\DN@{\let\@txtopt@\identity@#1}\fi\next@}\FN@\next@}
\def\check@therstyle#1#2{{\DN@{#1}\ifx\@txtopt@\identity@\else
		\DNii@##1\@therstyle{}\def\@therstyle{\DN@{#2}\nextii@}%
		\def\nextiii@##1##2\@therstyle{\expandafter\nextii@##1##2\@therstyle}%
    \expandafter\nextiii@\@txtopt@\@therstyle.\@therstyle\fi
	\expandafter}\next@}

\newif\if@theorem
\let\@therstyle\eat@
\def\@headtext@#1#2{{\disable@special\let\protect\noexpand
	\def\chapter{\protect\chapter@rh}\Protect@@\nofrills
	\edef\@tempa{\noexpand\@frills@\noexpand#1{#2}}\expandafter}\@tempa}
\let\AmSrighthead@\rightheadtext
\def\rightheadtext{\checkfrills@{\@headtext@\AmSrighthead@}}
\let\AmSlefthead@\leftheadtext
\def\leftheadtext{\checkfrills@{\@headtext@\AmSlefthead@}}
\def\@head@@#1#2#3#4#5{\if@theorem\else
	\@frills@{\csname\expandafter\eat@iv\string#4\endcsname}\@name{#1font@}\fi
  	\@name{#1\ST@LE}{\counter#3}{\ignorespaces#5\unskip}%
  \if@write\writeauxline{toc}{\eat@bs#1}{#2{\counter#3}{\ignorespaces#5}}\fi
	\if@theorem\else\expandafter#4\fi
	\ifx#4\endhead\ifx\@txtopt@\identity@\else
		\headmark{\@name{#1\ST@LE}{\counter#3}{\frills{}}}\fi\fi
	\ignorespaces}
\let\subsubheadfont@\empty
\ifx\undefined\endhead\Invalid@\endhead\fi
\def\@head@#1{\checkstar@{\checkfrills@{\checkbrack@{\@head@@#1}}}}
\def\@thm@@#1#2#3{%
	\@frills@{\csname\expandafter\eat@iv\string#3\endcsname}
  	{\@theoremtrue\check@therstyle{\@name{#1\ST@LE}}\frills
			{\counter#2}\@theoremfalse}%
	\expandafter\envir@stack\csname end\eat@bs#1\endcsname
	\@name{#1font@}\ignorespaces}
\let\proclaimfont@\empty
\def\@thm@#1{\checkstar@{\checkfrills@{\checkbrack@{\@thm@@#1}}}}
\def\@capt@@#1#2#3#4#5\endcaption{\bgroup
	\edef\@tempb{\global\footmarkcount@\the\footmarkcount@
    \global\@name{#2\N@M}\the\@name{#2\N@M}}%
	\def\shortcaption##1{\def\sh@rtt@xt####1{##1}}\let\sh@rtt@xt\identity@
	\DN@##1##2##3{\false@\fi\iftrue}%
	\ifx\@frills@\identity@\else\let\notempty\next@\fi
  #4{\@tempb\@name{#1\ST@LE}{\counter#2}}\@name{#1font@}#5\endcaption
  \if@write\writeauxline{#3}{\eat@bs#1}{{} \@name{#1\ST@LE}{\counter#2}%
    \if\notempty{#5}.\enspace\fi\sh@rtt@xt{#5}}\fi\egroup}
\def\@capt@#1{\checkstar@{\checkfrills@{\checkbrack@{\@capt@@#1}}}}
\let\captiontextfont@\empty

\def\newfont@def#1#2{\@ifundefined{#1font@}
	{\@xdefname{#1font@}{\@xname{.\expandafter\eat@iv\string#2font@}}}{}}
\def\newhead@#1#2#3#4{{%
	\gdef#1{\@therstyle\@therstyle\@head@{#1#2#3#4}}\newfont@def#1#4%
	\@ifundefined{#1\ST@LE}
    {\expandafter\gdef\csname\eat@bs#1\ST@LE\endcsname{\headstyle}}{}%
	\@ifundefined{#2}{\gdef#2{\headtocstyle}}{}%
  \@@addto\moretocdefs@{\\#1#1#4}}}
\outer\def\newhead#1{\checkbrack@{\expandafter\newhead@\expandafter
	#1\@txtopt@\headtocstyle}}
\outer\def\newtheorem#1#2#3#4{{%
	\gdef#2{\@thm@{#2#3#4}}\newfont@def#2\endproclaim%
	\@xdefname{\end\eat@bs#2}{\noexpand\revert@envir
		\@xname{\end\eat@bs#2}\noexpand#4}%
  \expandafter\gdef\csname\eat@bs#2\ST@LE\endcsname{\proclaimstyle{#1}}}}%
\outer\def\newcaption#1#2#3#4#5{{\let#2\relax
  \edef\@tempa{\gdef#2####1\csname end\eat@bs#2\endcsname}%
	\@tempa{\@capt@{#2#3{#4}#5}##1\endcaption}\newfont@def#2\endcaptiontext%
  \expandafter\gdef\csname\eat@bs#2\ST@LE\endcsname{\captionstyle{#1}}%
  \@@addto\moretocdefs@{\\#2#2\endcaption}\newtoc{#4}}}
{
\outer\gdef\newtoc#1{{%
	\expandafter\ifx\csname do@-#1\endcsname\relax
    \global\auxlinedef{#1}{\do@tocline{#1}}{}%
    \@@addto\tocsections@{\make@toc{#1}{}}\fi}}}

\toks@\expandafter{\itembox@}
\toks@@{{\let\therosteritem\identity@\let\rm\empty
  \edef\next@{\edef\noexpand\@lastmark{\therosteritem@}}\expandafter}\next@}
\edef\itembox@{\the\toks@@\the\toks@}
\def\firstitem@false{\let\iffirstitem@\iffalse
	\global\let\lastlabel\@lastlabel}

\def\rosteritemref#1{\hbox{\therosteritem{\@@@xref{#1}}}}
\def\local#1{\label@\@lastlabel{\lastlabel-i#1}}
\def\loccit#1{\rosteritemref{\lastlabel-i#1}}
\def\xRef@P@{\gdef\lastlabel}
\def\xRef#1{\@xref{#1}\protect\xRef@P@{#1}}

\def\iref@P@{\gdef\lastref}
\def\itemref#1#2{\rosteritemref{#1-i#2}\protect\iref@P@{#1}}
\def\iref#1{\@xref{#1}\itemref{#1}}
\def\ditto#1{\rosteritemref{\lastref-i#1}}

\def\eqref#1{\thetag{\@@@xref{#1}}}
\def\tagform@#1{\ifmmode\hbox{\rm\else\rom{\fi
	(\ignorespaces#1\unskip)\iftrue}\else}\fi}

\let\AmSfnote@\makefootnote@
\def\makefootnote@#1{{\let\footmarkform@\identity@
  \edef\next@{\edef\noexpand\@lastmark{#1}}\expandafter}\next@
  \AmSfnote@{#1}}

\def\clearpage{\ifnum\insertpenalties>0\line{}\fi\vfill\supereject}

\def\find@#1\in#2{\let\found@\false@
	\DNii@{\ifx\next\@nil\let\next\eat@\else\let\next\nextiv@\fi\next}%
	\edef\nextiii@{#1}\def\nextiv@##1,{%
    \edef\next{##1}\ifx\nextiii@\next\let\found@\true@\fi\FN@\nextii@}%
	\expandafter\nextiv@#2,\@nil}
\def\include@#1{{\ifx\@auxout\@subaux\DN@{\errmessage
		{\labelmesg@ Only one level of \string\include\space is supported}}%
	\else\edef\@tempb{#1}\@numopt@\clearpage
		\xdef\@include@{\relax\bgroup\ifx\@numopt@\identity@\clearpage
			\else\let\noexpand\@immediate\noexpand\empty\fi}%
	  \DN@##1 {\if\notempty{##1}\edef\@tempb{##1}\DN@####1\eat@ {}\fi\next@}%
  	\DNii@##1.{\edef\@tempa{##1}\DN@####1\eat@.{}\next@}%
		\expandafter\next@\@tempb\eat@{} \eat@{} %
  	\expandafter\nextii@\@tempb.\eat@.%
		\let\next@\empty
	  \if\expandafter\notempty\expandafter{\@tempa}%
		  \edef\nextii@{\write\@mainaux{%
  			\noexpand\string\noexpand\@input{\@tempa.aux}}}\nextii@
  		\@ifundefined\@includelist{\let\found@\true@}
				{\find@\@tempa\in\@includelist}%
			\if\found@\@ifundefined\@noincllist{\let\found@\false@}
				{\find@\@tempb\in\@noincllist}\else\let\found@\true@\fi
			\if\found@\@ifundefined{\cl@#1}{}{\DN@{\csname cl@#1\endcsname}}%
			\else\global\let\@auxout\@subaux
				\xdef\@auxname{\@tempa}\xdef\@inputname{\@tempb}%
	 	  	\@numopt@\immediate\openout\@auxout=\@tempa.aux
				\@numopt@\immediate\write\@auxout{\relax}%
 	  		\DN@{\@input\@inputname \@include@\closeaux@@\egroup\make@auxmain}\fi\fi\fi
  \expandafter}\next@}
\def\include{\checkstar@\include@}
\def\includeonly#1{\edef\@includelist{#1}}
\def\noinclude#1{\edef\@noincllist{#1}}

\newwrite\@mainaux
\newwrite\@subaux
\def\make@auxmain{\global\let\@auxout\@mainaux
  \xdef\@auxname{\jobname}\xdef\@inputname{\jobname}}
\begingroup
\catcode`\(\the\catcode`\{\catcode`\{=12
\catcode`\)\the\catcode`\}\catcode`\}=12
\gdef\closeaux@@((%
	\def\\##1(\@Waux(\global##1=\the##1))
	\edef\@tempa(\@Waux(%
		\gdef\expandafter\string\csname cl@\@auxname\endcsname{)%
		\\\pageno\\\footmarkcount@\closeaux@\@Waux(})\@immediate\closeout\@auxout)%
  \expandafter)\@tempa)
\endgroup
\let\closeaux@\empty
\let\@immediate\immediate
\def\@Waux{\@immediate\write\@auxout}
\def\readaux{%
	\W@{>>> \labelmesg@ Run this file twice to get x-references right}
	\global\everypar{}%
	{\def\\##1{\global\let##1\relax}%
		\def\/##1{\gdef##1{\preambule@cs##1}}%
		\disablepreambule@cs}%
	\make@auxmain\make@letter{\setboxz@h{\@input{\jobname.aux}}}\t@mpcat
  \immediate\openout\@auxout=\jobname.aux
	\immediate\write\@auxout{\relax}\global\let\end\endmain@}
\everypar{\global\everypar{}\readaux}
{\toks@\expandafter{\topmatter}
\global\edef\topmatter{\noexpand\readaux\the\toks@}}
\let\@@end@@\end

\def\endmain@{\clearpage\@immediate\closeout\@auxout
	\make@letter\def\newlabel##1##2##3{}\@input{\jobname.aux}%
 	\W@{>>> \labelmesg@ Run this file twice to get x-references right}%
 	\@@end@@}
\def\disablepreambule@cs{\\\disablepreambule@cs}
\def\preambule@cs#1{\warning@{Preamble command \string#1\space ignored}}

\def\proof{\checkfrills@{\checkbrack@{%
	\check@therstyle{\@frills@{\demo}{\frills{Proof}}{}}
		{\frills{}\envir@stack\endremark\envir@stack\enddemo}%
  \envir@stack\endproof\ignorespaces}}}
\def\endproof{\nofrillscheck{\frills@{\qed}\revert@envir\endproof\enddemo}}

\let\AmSref\ref
\let\AmSrefstyle\refstyle
\let\plaincite\cite
\def\citei@#1,{\citeii@#1\eat@,}
\def\citeii@#1\eat@{\@@@xref{#1}}
\def\cite#1{\protect\plaincite{\citei@#1\eat@,\unskip}}
\def\refstyle#1{\AmSrefstyle{#1}\uppercase{%
	\ifx#1A\relax \def\@ref@##1{\AmSref\xdef\@lastmark{##1}\key##1}%
  	\else\ifx#1C\relax \def\@ref@{\AmSref\no\counter\refno}%
		\else\def\@ref@{\AmSref}\fi\fi}}
\refstyle A
\newcounter\refno\null
\gdef\Refs{\checkstar@{\checkbrack@{\csname AmSRefs\endcsname
  \nofrills{\frills{References}%
  \if@write\writeauxline{toc}{vartocline}{\frills{References}}\fi}%
  \def\ref{\@ref@}\ignorespaces}}}
\let\ref\xref

\def\tocsections@{\make@toc{toc}{}}
\let\moretocdefs@\empty
\def\newtocline@#1#2#3{%
  \edef#1{\def\@xname{#2line}####1{\expandafter\noexpand
      \csname\expandafter\eat@iv\string#3\endcsname####1\noexpand#3}}%
  \@edefname{\no\eat@bs#1}{\let\@xname{#2line}\noexpand\eat@}%
	\@name{\no\eat@bs#1}}
\def\maketoc#1#2{\Err@{\Invalid@@\string\maketoc}}
\def\newtocline#1#2#3{\Err@{\Invalid@@\string\newtocline}}
\def\make@toc#1#2{\penaltyandskip@{-200}\aboveheadskip
	\if\notempty{#2}
		\centerline{\headfont@\ignorespaces#2\unskip}\nobreak
  	\vskip\belowheadskip \fi
	\@openin{#1}\@openout{#1}%
	\vskip\z@}
\def\contents{\readaux\checkfrills@{\checkbrack@{\@contents@}}}
\def\@contents@{\toc@{\frills{Contents}}\envir@stack\endcontents%
	\def\nopagenumbers{\let\page\eat@}\let\newtocline\newtocline@
  \def\tocline##1{\csname##1line\endcsname}%
  \edef\caption##1\endcaption{\expandafter\noexpand
    \csname head\endcsname##1\noexpand\endhead}%
	\ifmonograph@\def\vartoclineline{\Chapterline}%
		\else\def\vartoclineline##1{\sectionline{{} ##1}}\fi
  \let\\\newtocline@\moretocdefs@
	\ifx\@frills@\identity@\def\\##1##2##3{##1}\moretocdefs@
		\else\let\tocsections@\relax\fi
	\def\\{\unskip\space\ignorespaces}\let\maketoc\make@toc}
\def\endcontents{\tocsections@\vskip-\lastskip\revert@envir\endcontents
	\endtoc}

\@ifundefined\selectf@nt{\let\selectf@nt\identity@}{}
\def\textonlyfont@#1#2{%
	\expandafter\def\csname\eat@bs#1@P@\endcsname{%
  	\RIfM@\Err@{Use \string#1\space only in text}
		\else\edef\f@ntsh@pe{\string#1}\selectf@nt#2\fi}%
	\edef#1{\noexpand\protect\@xname{#1@P@}\empty}}
\tenpoint

\def\newshapeswitch#1#2{\gdef#1{\selectsh@pe#1#2}\PROTECT#1}
\def\shapeswitch#1#2#3{\expandafter\gdef\csname\eat@bs#1\string#2\endcsname{#3}}
\shapeswitch\rm\bf\bf  \shapeswitch\rm\tt\tt  \shapeswitch\rm\smc\smc
\newshapeswitch\em\it
\shapeswitch\em\it\rm  \shapeswitch\em\sl\rm
\def\selectsh@pe#1#2{\relax\@ifundefined{#1\f@ntsh@pe}{#2}
	{\csname\eat@bs#1\f@ntsh@pe\endcsname}}

\def\@itcorr@{\leavevmode\skip@\lastskip\unskip\/%
  \ifdim\skip@=\z@\else\hskip\skip@\fi}
\def\rom@P@#1{\@itcorr@{\selectsh@pe\rm\rm#1}}
\def\rom{\protect\rom@P@}
\def\Rom@P@#1{\@itcorr@{\rm#1}}
\def\Rom{\protect\Rom@P@}
{\catcode`\-=11
\gdef\wr@index#1{\@writeaux\tf@-idx{\string\indexentry{#1}{\thepage}}}
\gdef\wr@glossary#1{\@writeaux\tf@-glo{\string\glossaryentry{#1}{\thepage}}}}
\def\emph{\@itcorr@\bgroup\em\checkstar@\emph@}
\let\index\wr@index
\let\glossary\wr@glossary
\def\emph@#1{\@numopt@{\wr@index{#1}}\ignorespaces#1\unskip\egroup
  \DN@{\DN@{}\ifx\next.\else\ifx\next,\else\DN@{\/}\fi\fi\next@}\FN@\next@}
\def\makequoteactive{\catcode`\"\active}
{\makequoteactive\gdef"#1"{\emph{#1}}}
\def\MakeIndex{\@openout{idx}}
\def\MakeGlossary{\@openout{glo}}

\def\endofpar#1{\ifmmode\ifinner\endofpar@{#1}\else\eqno{#1}\fi
	\else\leavevmode\endofpar@{#1}\fi}
\def\endofpar@#1{\unskip\penalty\z@\null\hfil\hbox{#1}\hfilneg\penalty\@M}

\@@addto\disablepreambule@cs{%
	\\\readaux
	\/\Monograph
	\/\MakeIndex
	\/\MakeGlossary
}

\catcode`\@\labelloaded@


\def\headstyle#1#2{\numberline{#1.\enspace}#2}
\def\headtocstyle#1#2{\numberline{#1.}\space #2}

\def\specialtocstyle#1#2{#2}
\newcounter\section\null
\newcounter\subsection\section
\newcounter\subsubsection\subsection
\newhead\specialsection[\specialtocstyle]\null\endspecialhead
\newhead\section\section\endhead
\newhead\subsection\subsection\endsubhead
\newhead\subsubsection\subsubsection\endsubsubhead
\def\firstappendix{\global\sectionno=0 %
  \counterstyle\section{\Alphnum\sectionno}%
	\global\let\firstappendix\empty}

\def\appendixtocstyle#1#2{\space\numberline{Appendix #1.\enspace}#2}
\newhead\appendix[\appendixtocstyle]\section\endhead

\let\endAmSdef\enddefinition
\def\proclaimstyle#1#2{\numberline{#2.\enspace}\frills{#1}}
\copycounter\thm\subsubsection
\theorem\thm\endproclaim
\proposition\thm\endproclaim
\lemma\thm\endproclaim
\corollary\thm\endproclaim
\definition\thm\endAmSdef
\example\thm\endAmSdef

\def\captionstyle#1#2{\frills{\numberline{#1 #2}}}
\newcounter\figure\null
\newcounter\table\null
\newcaption{Figure}\figure\figure{lof}\botcaption
\newcaption{Table}\table\table{lot}\topcaption

\copycounter\equation\subsubsection

\immediate\write10{This is DEGT.DEF by A.Degtyarev <22 October 1995>}
\chardef\tempcat\the\catcode`\@
\catcode`\@=11
%
%
%
\let\ge\geqslant
\let\le\leqslant
%
\def\C{{\Bbb C}}
\def\R{{\Bbb R}}
\def\Z{{\Bbb Z}}

%
%
%
\def\Cp#1{\C{\operator@font p}^{#1}}
\def\Rp#1{\R{\operator@font p}^{#1}}
%
\def\Hom{\qopname@{Hom}}
\def\Ext{\qopname@{Ext}}
\def\Tors{\qopname@{Tors}}

\def\Im{\qopname@{Im}}			
\def\Re{\qopname@{Re}}     %
\def\Ker{\qopname@{Ker}}
\def\Coker{\qopname@{Coker}}
\def\Int{\qopname@{Int}}
\def\Cl{\qopname@{Cl}}
\def\Fr{\qopname@{Fr}}
\def\Fix{\qopname@{Fix}}
\def\tr{\qopname@{tr}}
\def\inj{\qopname@{in}}
\def\id{\qopname@{id}}
\def\pr{\qopname@{pr}}
\def\rel{\qopname@{rel}}
\def\pt{{\operator@font{pt}}}
\def\const{{\operator@font{const}}}
\def\codim{\qopname@{codim}}
\def\cdim{\qopname@{dim_{\C}}}
\def\rdim{\qopname@{dim_{\R}}}
\def\conj{\qopname@{conj}}
\def\rank{\qopname@{rk}}
\def\sign{\qopname@{sign}}
\def\gcd{\qopname@{g.c.d.}}

\def\set<#1|#2>{\bigl\{#1\bigm|#2\bigr\}}
\def\emptyset{\varnothing}
%

%

%
%
\def\preprint#1{\hrule height0pt depth0pt\kern-24pt%
  \hbox to\hsize{#1}\kern24pt}
\def\today{\ifcase\month\or January\or February\or March\or
  April\or May\or June\or July\or August\or September\or October\or
  November\or December\fi \space\number\day, \number\year}
\def\n@te#1#2{\leavevmode\vadjust{%
 {\setbox\z@\hbox to\z@{\strut\eightpoint#1}%
  \setbox\z@\hbox{\raise\dp\strutbox\box\z@}\ht\z@=\z@\dp\z@=\z@%
  #2\box\z@}}}
\def\leftnote#1{\n@te{\hss#1\quad}{}}
\def\rightnote#1{\n@te{\quad\kern-\leftskip#1\hss}{\moveright\hsize}}
\def\?{\FN@\qumark}
\def\qumark{\ifx\next"\DN@"##1"{\leftnote{\rm##1}}\else
 \DN@{\leftnote{\rm??}}\fi{\rm??}\next@}
\def\centerpage{\dimen@=6.5truein \advance\dimen@-\hsize\hoffset.5\dimen@}
\ifnum\mag>1000 \centerpage\fi
%
\def\nologo{\let\logo@\relax}
%
%
\def\linefill@@{\mkern-2mu\mathord-\mkern-2mu}
\def\Linefill@@{\mkern-2mu\mathord=\mkern-2mu}
\expandafter\ifx\csname msa@group\endcsname\relax 
 \expandafter\ifx\csname msafam\endcsname\relax\let\next\relax 
 \else
  \edef\next{\hexnumber@\msafam}
  \font@\eightmsa=msam8
  \font@\sixmsa=msam6
  \addto\eightpoint{\textfont\msafam\eightmsa \scriptfont\msafam\sixmsa}
  \addto\tenpoint{\textfont\msafam\tenmsa \scriptfont\msafam\sevenmsa}
 \fi\else
   \edef\next{\ifcase\msa@group0 \or1 \or2 \or3 \or4 \or5 \or6 \or7
      \or8 \or9 \or A \or B \or C \or D \or F \fi}
\fi
\ifx\next\relax\let\setd@sh\relax\else
 \mathchardef\d@shright"3\next4B
 \mathchardef\d@shleft"3\next4C
 \def\setd@sh{{\ifnum\mindasharrowwidth=0 \else 
  \setboxz@h{$\m@th\dabar@$}\dimen@\wdz@\multiply\dimen@\mindasharrowwidth
  \loop\ifdim\dimen@<\bigaw@\advance\dimen@\wdz@\repeat
  \global\bigaw@\dimen@\fi}}
\fi
\newcount\mindasharrowwidth
\mindasharrowwidth=3
\newif\ifdouble@rrow
\newif\ifleft@rrow
\newif\ifd@sh
\def\arrowhead{{\relax\ifd@sh\ifleft@rrow\d@shleft\else\d@shright\fi
 \else\ifdouble@rrow\ifleft@rrow\Leftarrow\else\Rightarrow\fi
  \else\ifleft@rrow\leftarrow\else\rightarrow\fi\fi\fi}}
\def\plaintail{{\relb@r}}
\def\arrowtail{{\relax\ifleft@rrow\mapstochar\relb@r
 \else\relb@r\mskip-2.5mu\mapstochar\mskip 2.5mu\fi}}
\def\cl@p#1{\mathrel{\hbox to\z@{\hss$\m@th#1$\hss}}}
\def\hook{{\relax\ifleft@rrow\cl@p{\lhook\joinrel}\relb@r\else
 \relb@r\cl@p{\mskip-1.5mu\rhook}\fi}}
\def\doublehead{{\relax\ifd@sh
  \ifleft@rrow\cl@p{\d@shleft\joinrel}\d@shleft
   \else\d@shright\cl@p{\joinrel\d@shright}\fi
 \else\arrowhead\mkern-14mu\arrowhead\fi}}
\def\hexpand(#1,#2){\global\def\lsp@n{\dimen@#1\kern-\dimen@}%
 \global\def\rsp@n{\dimen@#2\kern-\dimen@}%
 \global\def\hsp@n{\global\bigaw@\z@\mindasharrowwidth\z@}}
\def\hexpand(#1,#2){\global\def\lsp@n{\dimen@#1\global\let\lsp@n\relax
 \ifCD@\kern-\dimen@\def\rsp@n{\dimen@#2\kern-\dimen@}
  \def\hsp@n{\global\bigaw@\z@\mindasharrowwidth\z@}
 \else\global\advance\bigaw@\dimen@\advance\minaw@#2\fi}}
\let\lsp@n\relax
\let\rsp@n\relax
\let\hsp@n\relax
\def\CDline@@{\lsp@n\ifCD@\enskip$\global\bigaw@\minCDaw@
  \else\global\bigaw@\minaw@\mathrel{\fi
 \setboxz@h{$\m@th\ssize\;{\s@p}\;\;$}
 \setbox@ne\hbox{$\m@th\ssize\;{\s@b}\;\;$}
 \setbox\tw@\hbox{$\m@th\s@b$}
 \ifdim\wdz@>\bigaw@\global\bigaw@\wdz@\fi
 \ifdim\wd@ne>\bigaw@\global\bigaw@\wd@ne\fi
 \hsp@n\ifd@sh\setd@sh\else\ifCD@\setd@sh\fi\fi\hskip.5\bigaw@ plus.5fil
 \hbox to\z@{\hss$\m@th\mathop{\vphantom=}\limits^{\s@p}
   \ifdim\wd\tw@>\z@ _{\s@b}\fi$\hss}
 \hskip-.5\bigaw@ plus-.5fil\rlap{\left@rrowtrue$\m@th\l@ft$}
 \ifd@sh\else\global\advance\bigaw@-6\ex@\kern 3\ex@\fi
 \cleaders\hbox{$\m@th\fill@$}\hskip\bigaw@ plus1fil
 \ifd@sh\else\kern 3\ex@\fi\llap{\left@rrowfalse$\m@th\r@ght$}
 \ifCD@$\enskip\else}\fi\rsp@n\d@shfalse\global\let\hsp@n\relax
 \global\let\lsp@n\relax\global\let\rsp@n\relax\ampersand@}
\def\CDline@@{\ifCD@$\global\bigaw@\minCDaw@
  \else\global\bigaw@\minaw@\mathrel{\fi\lsp@n
 \setboxz@h{$\m@th\ssize\;{\s@p}\;\;$}
 \setbox@ne\hbox{$\m@th\ssize\;{\s@b}\;\;$}
 \setbox\tw@\hbox{$\m@th\s@b$}
 \ifdim\wdz@>\bigaw@\global\bigaw@\wdz@\fi
 \ifdim\wd@ne>\bigaw@\global\bigaw@\wd@ne\fi\hsp@n
 \ifCD@\enskip\setd@sh\else\ifd@sh\setd@sh\fi\fi\hskip.5\bigaw@ plus.5fil
 \hbox to\z@{\hss$\m@th\mathop{\vphantom=}\limits^{\s@p}
   \ifdim\wd\tw@>\z@ _{\s@b}\fi$\hss}
 \hskip-.5\bigaw@ plus-.5fil\rlap{\left@rrowtrue$\m@th\l@ft$}
 \ifd@sh\else\global\advance\bigaw@-6\ex@\kern 3\ex@\fi
 \cleaders\hbox{$\m@th\fill@$}\hskip\bigaw@ plus1fil
 \ifd@sh\else\kern 3\ex@\fi\llap{\left@rrowfalse$\m@th\r@ght$}
 \ifCD@\enskip\rsp@n$\else}\fi\d@shfalse\ampersand@}
\def\ht@st#1{\FN@\ht@@st#1}
\def\ht@@st{\ifx\next<\DN@{\arrowhead\eat@}\else
 \ifx\next>\DN@{\arrowhead\eat@}\else
 \ifx\next-\DN@{\plaintail\eat@}\else
 \ifx\next=\DN@{\plaintail\eat@}\else
 \ifx\next:\DN@{\plaintail\eat@}\else
 \ifx\next|\DN@{\arrowtail\eat@}\else
 \ifx\next(\DN@{\hook\eat@}\else
 \ifx\next)\DN@{\hook\eat@}\else
 \let\next@\relax\fi\fi\fi\fi\fi\fi\fi\fi\next@}
\def\at@@st{\ifx\next-\let\next\CDline@
 \else\ifx\next:\let\next\CDdash@\else\let\next\CDLine@\fi\fi\next}
\newcount\spc@unt
\def\CDline@-#1-#2-#3{\def\s@p{#1}\def\s@b{#2}\def\r@ght{\ht@st#3}%
 \double@rrowfalse\let\relb@r\relbar\let\fill@\linefill@@\CDline@@}
\def\CDLine@=#1=#2=#3{\def\s@p{#1}\def\s@b{#2}\def\r@ght{\ht@st#3}%
 \double@rrowtrue\let\relb@r\Relbar\let\fill@\Linefill@@\CDline@@}
\def\CDdash@:#1:#2:#3{\def\s@p{#1}\def\s@b{#2}\def\r@ght{\ht@st#3}%
 \d@shtrue\def\relb@r{\mathrel\dabar@}\let\fill@\dabar@\CDline@@}
\def\CDhor@#1{\afterassignment\CDhor@@\global\spc@unt0#1}
\def\CDhor@@#1{\ampersand@\ifCD@\multispan{\the\spc@unt}\fi%
 \def\l@ft{\ht@st#1}\FN@\at@@st}
\atdef@:{\relax\ifmmode\CDhor@\else\leavevmode\null:\fi}
\atdef@>#1>#2>{\CDhor@--{#1}-{#2}->}
\atdef@<#1<#2<{\CDhor@<-{#1}-{#2}--}
\atdef@={\CDhor@=====}
\atdef@[#1>#2>{\CDhor@(-{#1}-{#2}->}
\atdef@]#1<#2<{\CDhor@<-{#1}-{#2}-)}
%
%
\def\pr@fix{\Big}
\def\short{\def\pr@fix{\big}}
\def\CDarrowsize#1{\def\pr@fix{#1}}
\def\ampersand@@{\ifx\next$\let\next@\relax\else
 \DN@{\ampersand@\ampersand@}\fi\next@}
\let\vert@ver\relax
\let\vert@nder\relax
\def\@rrow{\ifdouble@rrow\mathchar"377\else\mathchar"33F\fi}
\let\@rrowup\@rrow
\let\@rrowdown\@rrow
\def\set@up{\def\@rrowup{\ifdouble@rrow\mathchar"37E\else\mathchar"378\fi}}
\def\set@down{\def\@rrowdown{\ifdouble@rrow\mathchar"37F\else\mathchar"379\fi}}
\def\vertd@sh{\def\@rrow{\vbox{\kern.5\ex@\hbox to\wdz@
 {\hss\vrule width.4\ex@ height3.5\ex@\hss}\kern.5\ex@}}}
\def\vertd@sh{\edef\@rrow{\vbox{\kern.5\ex@\hbox{$\m@th\@rrow$}\kern.5\ex@}}}
\atdef@+{\FN@\CDvert@}
\atdef@ V#1V#2V{\CDvert@@<#1|\set@down|#2>}
\atdef@ A#1A#2A{\CDvert@@<#1|\set@up|#2>}
\atdef@|{\double@rrowtrue\FN@\CDvert@}
\atdef@\vert{\double@rrowtrue\FN@\CDvert@}
\def\CDvert@{\ifx\next<\let\next@\CDvert@@\else\DN@{\CDvert@@<||>}\fi\next@}
\def\CDvert@@<#1|#2|#3>{\setboxz@h{$\pr@fix.$}\nousp@n\nodsp@n
 \advance\dimen@ii\dp\z@\advance\dimen@\ht\z@\advance\dimen@\dimen@ii
 \llap{$\m@th\vcenter{\hbox{$\ssize#1$}}$}
 \boxz@\setboxz@h{\vbox to\dimen@{\offinterlineskip
  \let\up\set@up\let\down\set@down\let\double\double@rrowtrue#2
  \setboxz@h{$\m@th\@rrowdown$}\hbox to\wdz@{\hss\smash{\vert@ver}\hss}
  \hbox{$\m@th\@rrowup$}\kern-3\ex@\cleaders\hbox{$\m@th\@rrow$}\vfil
  \kern-3\ex@\copy\z@\hbox to\wdz@{\hss\smash{\vert@nder}\hss}}}
 \dp\z@-\dimen@ii\ht\z@\dimen@ii\raise-\dimen@ii\boxz@
 \rlap{$\m@th\vcenter{\hbox{$\ssize#3$}}$}\FN@\ampersand@@}
\def\upmapsto{\def\vert@ver{\vrule height.4\ex@ width2.4\ex@}}
\def\downmapsto{\def\vert@nder{\vrule depth.4\ex@ width2.4\ex@}}
%
\def\h@@k#1#2#3#4{\def#1{\raise#2\ex@\hbox{$\m@th\scriptscriptstyle\mkern#3mu#4$}}}
\def\ulhook{\h@@k\vert@ver{-2}{8}\cap}
\def\urhook{\h@@k\vert@ver{-2}{-7.1}\cap}
\def\dlhook{\h@@k\vert@nder{-1}{8}\cup}
\def\drhook{\h@@k\vert@nder{-1}{-7.1}\cup}
\def\nousp@n{\dimen@\z@}
\def\nodsp@n{\dimen@ii\z@}
\def\vexpand(#1,#2){\ifCD@\dimen@#1\dimen@ii#2\let\nousp@n\relax\let\nodsp@n\relax\fi}
\def\expand{\ifCD@\dimen@ii 2.5\ex@\let\nodsp@n\relax\fi}
\def\Expand{\ifCD@\dimen@ii 21.5\ex@\let\nodsp@n\relax\fi}
\def\expandup{\ifCD@\dimen@ 2.5\ex@\let\nousp@n\relax\fi}
\def\Expandup{\ifCD@\dimen@ 21.5\ex@\let\nousp@n\relax\fi}
%
%
%
\expandafter\ifx\csname minCDaw@\endcsname\relax
 \let\minCDaw@\minCDarrowwidth\fi
\expandafter\ifx\csname eat@\endcsname\relax\def\eat@#1{}\fi
\expandafter\ifx\csname ssize\endcsname\relax\let\ssize\scriptstyle\fi
\expandafter\ifx\csname operator@font\endcsname\relax
 \def\operator@font{\roman}\fi
\expandafter\ifx\csname eightpoint\endcsname\relax
 \let\eightpoint\small\fi
\catcode`\@\tempcat

\def\mnote#1{\rightnote{\rm#1}}
\Protect\mnote

\def\+{\mathbin{\scriptstyle\sqcup}}
\let\<\langle
\let\>\rangle

\Protect\cong

\def\tps-#1{\vbox{\hbox{\rm S-#1}}}
\def\tph#1-#2{\vbox{\hbox{\rm H#1-#2}}}
\def\tpc#1-#2{\vbox{\hbox{\rm C#1-#2}}}
\def\tpi#1-#2{\vbox{\hbox{\rm I#1-#2}}}

\def\er{E_{\R}\futurelet\nexts\getsup}
\def\xr{X_{\R}\futurelet\nexts\getsup}
\Protect\er
\Protect\xr

\def\c{\futurelet\nexts\conjs}
\Protect\c
\def\index#1{^{\botsmash{(#1)}}}
\def\conjs{\ifx\nexts1\def\nexts##1{t\index1}\else
 \ifx\nexts2\def\nexts##1{t\index2}\else
 \def\nexts{\mathop{\roman{conj}}}\fi\fi\nexts}
\def\eri{\er\index i}
\def\xri{\xr\index i}
\Protect\eri
\Protect\xri

\def\getsup{\ifx\nexts1\def\nexts##1{\index1}\else
 \ifx\nexts2\def\nexts##1{\index2}\else\let\nexts\relax\fi\fi\nexts}

\def\subindex#1{_{\topsmash{(#1)}}}
\def\subtwo{\subindex2}
\let\2\subtwo

\def\I{\roman{I}}

\def\II{\roman{II}}

\Protect\Cal

\def\discr{\operatorname{discr}}

\def\Pic{\operatorname{Pic}}

\def\tP{\smash{\tilde P}}
\def\O{\Cal O}

\def\t{\tau}
\def\cy{c}
\def\yrr{Y_\R\futurelet\nexts\getsup}
\def\crr{C_\R\futurelet\nexts\getsup}
\def\yrri{Y_\R\index i}
\def\crri{C_\R\index i}
\def\coord(#1:#2){(#1\,{:}\,#2)}

\readaux
\preprint{\hss\eightrm Submitted to \eightit Math. Ann.}

\topmatter
\title
Empty real Enriques surfaces and Enriques-Einstein-Hitchin 4-manifolds
\endtitle
\rightheadtext{Empty Enriques surfaces} 

\author
Alexander Degtyarev and Viatcheslav Kharlamov
\endauthor

\address
Steklov Mathematical Institute, St.~Petersburg branch, Russia
\newline\indent
and Bilkent University, Turkey
\endaddress

\email
degt\,\@\,fen.bilkent.edu.tr, degt\,\@\,pdmi.ras.ru
\endemail

\address
Institut de Recherche Math\'ematique Avanc\'ee,\newline\indent
Universit\'e Louis Pasteur et CNRS, France
\endaddress

\email
kharlam\,\@\,math.u-strasbg.fr
\endemail




\endtopmatter

\section{Introduction}\label{intro}
N.~Hitchin~\cite{Htchn} proved that the Euler
characteristic~$\chi(E)$ and
signature~$\sigma(E)$ of a compact orientable $4$-dimensional
Einstein manifold~$E$ satisfy the inequality
$|\sigma(E)|\le\frac23\chi(E)$, the equality holding only if
either $E$ is flat or the universal covering~$X$ of~$E$ is a
$K3$-surface and $\pi_1(E)=1$, $\Z/2$, or $\Z/2\times\Z/2$. In the
latter cases, $E$ is a $K3$-surface if $\pi_1=1$, an Enriques surface
if $\pi_1=\Z/2$, or the quotient of an Enriques surface by a free
antiholomorphic involution if $\pi_1=\Z/2\times\Z/2$. It is the
Einstein manifolds of the last type that we call
\emph{Enriques-Einstein-Hitchin varieties}. The varieties of the
other three extremal types (flat, $\pi_1=1$, and $\pi_1=\Z/2$) are
known to form connected families: two varieties of the same type can
be deformed continuously into each other. To our knowledge, the
number of connected components of the moduli space of
Enriques-Einstein-Hitchin varieties was not known. In this paper we
give the answer: we prove that {\sl their moduli space is connected.}

As is known (modulo Calabi-Yau theorem this statement is also
contained in~\cite{Htchn}), the universal covering~$X$ of an
Enriques-Einstein-Hitchin manifold~$E$ carries a canonical complex
structure, so that $X$ is a $K3$-surface, one nontrivial element of
$\pi_1(E)=\Z/2\times\Z/2$ acts on~$X$ holomorphically, and the two
others, anti-holomorphically. This correspondence establishes a
homotopy equivalence between the moduli space of
Enriques-Einstein-Hitchin varieties and that of Enriques surfaces
with free anti-holomorphic involution (cf.~\cite{Itoh}).  An Enriques
surface with a free anti-holomorphic involution is, by definition, an
empty real Enriques surface, and the connectedness of the moduli
space of Enriques-Einstein-Hitchin varieties follows from the main
result of the present paper:

\theorem[\subsection{Main Theorem}]\label{main}
All empty real Enriques surfaces \rom(or, equivalently, compact
Einstein $4$-manifolds with $\vert\sigma\vert = \frac23 \chi$ and
$\pi_1=\Z/2\times\Z/2$\rom) are of the same deformation type.
\endtheorem

Originally this result was obtained as part of the solution of a more
general
problem:  we enumerated the connected components of the moduli space
of all (not only empty) real Enriques surfaces. For this purpose we
developed two different approaches:  one is based on what we call
Donaldson-Hitchin trick and reduces the task to a geometrical study
of real rational surfaces, the other one uses an explicit description
of the moduli space and requires an arithmetical study of
$(\Z/2\times\Z/2)$-actions in the intersection lattice of a
$K3$-surface, see~\cite{cras}. The complete solution was obtained in
collaboration with I.~Itenberg; it will appear in~\cite{DIK}.

In this paper we present a different proof; it is based on a
systematic study of real elliptic pencils and gives explicit models
of all empty real Enriques surfaces.  It has many similarities with
(and relies upon) Horikawa's proof of connectedness of the moduli
space of complex Enriques surfaces (see~\cite{Hor}) and with Cosec's
and Dolgachev's study of complex elliptic pencils on Enriques
surfaces (see~\cite{DC} and~\cite{DinInv84}).

\section{Preliminaries}\label{pre}

In this section we fix main definitions, introduce principal
notation and recall some known basic results (most of them are
found, e.g., in ~\cite{BPV} and~\cite{DC}).

\subsection{$K3$- and Enriques surfaces}\label{2.1}
A \emph{$K3$-surface} is a compact complex analytic surface~$X$ with
$\pi_1(X)=1$ and $c_1(X)=0$.
An \emph{Enriques surface} is the quotient of a K3-surface by a fixed
point free holomorphic involution. A \emph{real Enriques surface} is
an Enriques surface equipped with a \emph{real structure}, i.e.,
anti-holomorphic involution.  The real structure lifts to
the covering $K3$-surface (for a fixed point free involution a proof
is given in~\cite{Htchn}; in others cases it is straightforward).
Thus, real Enriques surfaces are in a one-to-one correspondence with
$K3$-surfaces equipped with an action of a group
$G\cong\Z/2\times\Z/2$ generated by two commuting
anti-holomorphic involutions with fixed point free composition.
(Note that the isomorphism
$G\cong\Z/2\times\Z/2$ is not fixed, cf. Remark in~\ref{last}.)

For an Enriques surface~$E$ the map $\Pic E\to H_2(E;\Z)$,
$D\mapsto[D]$, is an isomorphism. The nontrivial element of
$\Tors H_2(E;\Z)=\Z/2$ is equal to $[K_E]=w_2(E)$, where $K_E\in\Pic E$ is
the canonical class.  Denote
$L=H_2(E;\Z)/\Tors$; it is a lattice isomorphic to $E_8\oplus U$.
Here, $E_8$ is the negative lattice generated by the root system of
the same name and $U$
is a hyperbolic plane.  In $U$ there are two pairs of generators
$x_1$, $x_2$ with $x_1^2=x_2^2=0$, $x_1 x_2=1$; we call such
generators \emph{standard}.

Throughout the paper we denote by~$E$ a fixed Enriques surface,
by~$X$, its universal covering, and by $\c\:E\to E$, the real
structure on~$E$ (if $E$ is real).
The eigenlattices
of~$\c_*\:L\to L$ are denoted by~$L^\pm$. For an element $y\in L$ we
use~$\hat y$ to denote one of the two lifts of~$y$ to $H_2(E;\Z)$.
Given $y\in L^-$ or $x\in H_2(E;\Z)$ with $(x\bmod\Tors)\in L^-$, we
let $\delta(y)=(1+\c_*)\hat y$ and $\delta(x)=(1+\c_*)x$.  Clearly,
$\delta$ can take only two values:  $0$ and~$w_2(E)$.

As usual, we denote by~$|D|$ the linear system generated by a
divisor~$D$. In the case of Enriques surface the Riemann-Roch theorem
takes the form
$$
\tsize
\dim|D|-\dim H^1(E;\O_E(D))+\dim|K_E-D|=\frac12D^2.
$$
Since, in addition, $K_E$ is not effective (as $2K_E=0$), this
implies that for any element $x\in H_2(E;\Z)$ with $x^2\ge0$
either~$x$ or~$-x$ (but not both) is realized by an effective
divisor.

A nonsingular rational curve~$R$ on an Enriques surface~$E$ (or
$K3$-surface~$X$) is called a \emph{nodal curve}. As follows from the
adjunction formula, $R^2=-2$. Hence, nodal curves are determined by
their classes in $H_2(E;\Z)$ (respectively, $H_2(X;\Z)$): two nodal
curves of the same class coincide. In the case of Enriques surface
$R$ is, in fact, determined by its class in~$L$; if $E$ is real
and $([R]\bmod\Tors)\in L^-$, then $\delta(R)=0$. An Enriques or
$K3$-surface is called \emph{nodal} (respectively, \emph{unnodal}) if
it has (respectively, does not have) a nodal curve.
Note that `most' Enriques surfaces are unnodal, cf. the proof
of~\ref{gen}.

\subsection{Elliptic pencils}\label{2.2}
An \emph{elliptic pencil} on an algebraic surface is a linear system
of dimension one whose generic member is an irreducible smooth elliptic
curve. Recall (see~\cite{DC} or \cite{BPV}, Lemma 17.1) that any
elliptic pencil on an Enriques surface~$E$ is base point free and,
thus, defines an elliptic fibration of $E\to\Cp1$ with exactly
two multiple fibers~$P$, $P'$, so that $\dim|P|=\dim|P'|=0$ and
$P-P'=K_E$. The pencil is then the linear system $|2P|=|2P'|$.

Let $|2P|$ be an elliptic pencil on~$E$ and $P$, $P'$ its multiple
fibers. Their pull-backs~$\tP$, $\tP'$ in~$X$,
which are irreducible smooth elliptic curves, are members of an
elliptic pencil on~$X$, which we denote by~$|\tP|$. The corresponding
fibration is induced from $|2P|$ by the double covering of its base
branched at the points corresponding to~$P$ and~$P'$. We will
call~$|\tP|$ the \emph{pull-back of~$|2P|$} and~$|2P|$, the
\emph{projection of~$|\tP|$}.

The pull-back~$|\tP|$ is $\tau$-equivariant.
Conversely, if $|\tP|$ is a $\tau$-invariant elliptic pencil on~$X$,
the projections of its members to~$E$ have self-intersection~$0$ and
are homologous and, hence, linear equivalent to each other; thus,
they form an elliptic pencil. (In particular, $\tau$
always acts nontrivially on the base of an equivariant pencil, as
otherwise one would obtain an elliptic pencil in~$E$ without multiple
fibers.)

An elliptic pencil on a real Enriques surface is called \emph{real} if
it is $\c$-invariant. Clearly, the pull-back of a real elliptic pencil is
real.

\proposition\label{reality}
An elliptic pencil $|2P|$ is real if and only if
$([P]\bmod\Tors)\in L^-$. If this is the case,
the two multiple fibers, $P$ and $P'$, of the pencil are real if
$\delta([P])=0$ and conjugate to each other if $\delta([P])=w_2(E)$.
\endproposition

\proof
Since $\c$ is antiholomorphic, it transforms holomorphic curves to
holomorphic curves reversing their complex orientation (and their
complex structure). Thus, the first statement follows from
$\Pic(E)=H_2(E;\Z)$.  The second statement follows from $P-P'=K_E$,
where $P$ and $P'$ are the multiple fibers.
\endproof

\proposition\label{real.part}
The involution induced by~$\conj$ on the base of a real elliptic
pencil on a real Enriques surfaces has nonempty real part.
\endproposition

\proof
Let $|2P|$ be the pencil, $B$ its base, and $|\tP|$ the pull-back of
$|2P|$ in~$X$. The base~$\tilde B$ of~$|\tP|$ is the double covering
of~$B$ branched at the two points corresponding to the two multiple
fibers. Since $|\tP|$ is also a real elliptic pencil, the involution
induced by~$\c$ on~$B$ must lift to an involution on~$\tilde B$. On
the other hand, any real structure on~$B$ with $B_\R=\varnothing$
lifts to an order~$4$ automorphism of~$\tilde B$.
\endproof

\section{Models}\label{s3}

\subsection{Nonspecial Horikawa models}\label{3.1} Let $X$ be the
$K3$-surface obtained as the minimal resolution of the double
covering of $Y=\Cp1\times\Cp1$ branched over a reduced
bi-degree~$(4,4)$ curve~$C\subset Y$ with at worst simple
singularities. Denote by $s\:Y\to Y$ the Cartesian product of the
nontrivial involutions $\coord(u:v)\mapsto\coord(-u:v)$ of the
factors. Up to isomorphism $s$ is the only holomorphic involution
on~$Y$ with isolated fixed points. If $C$ is $s$-symmetric, $s$ lifts
to two different involutions on~$X$, commuting with the deck
translation~$d$ of $X\to Y$. If, besides, $C$ contains no fixed
points of~$s$, exactly one of these two involutions, which we denote
by~$\t$, is fixed point free (see, e.g.,~\cite{Hor} or~\cite{BPV}),
and, hence, the orbit space $E=X/\t$ is an Enriques surface. The two
rulings of $Y$ define two $\tau$-invariant elliptic pencils on $X$,
which project to two elliptic pencils on $E$. We call these pencils
{\it basic}; their multiple fibers correspond to the generatrices
of~$Y$ through the fixed points of~$s$.

Suppose that $Y$ is equipped with a real structure $\cy$ commuting
with~$s$ and $C$ is a real curve. Then $s\circ\cy$ is another real
structure on~$Y$ and~$C$.  We denote the real point sets of these
structures by $\yrri$ and $\crri$, $i=1,2$ ($i=1$ corresponding to
$\cy$) and call them the {\it halves\/} of~$Y$ and~$C$.  The
involutions~$\cy$ and~$s\circ\cy$ lift to four different commuting
real structures ($\c1$, $\c2=\t\circ\c1$, $d\circ\c1$, and
$d\circ\c2$) on~$X$, which, in turn, descend to two real structures
on~$E$, called the {\it expositions\/} of~$E$.  A choice of an
exposition is determined by a choice of one of the two liftings
$\c1$, $d\circ\c1$ of~$\cy$ to~$X$.

The involutions $s$ and $\cy$ generate a $\Z/2\times\Z/2$-action
on~$Y$, so that one nontrivial element acts holomorphically and has
isolated fixed points and two other nontrivial elements act
anti-holomorphically.

\proposition\label{actions}
Up to isomorphism, there are five such actions. Four of them are
\emph{decomposable}, i.e., split into product of actions on the
factors\rom:
\roster
\item\local1
$\yrr1$, $\yrr2$ are homeomorphic to $S^1\times S^1$\!,
each ruling has two invariant fibers\rom;
\item\local2
$\yrr1$ is $S^1\times S^1$, $\yrr2=\varnothing$, and one ruling
has two invariant fibers\rom;
\item\local3
$\yrr1$ is $S^1\times S^1$, $\yrr2=\varnothing$, and there is no
invariant fibers\rom;
\item\local4
$\yrr1=\yrr2=\varnothing$ and there is no invariant fibers\rom;
\endroster
and one is\/ \emph{indecomposable}:
\roster
\item[5]\local5
$\yrr1$ and $\yrr2$ are homeomorphic to $ S^2$, and $s$ has two real
fixed points.
\endroster
\endproposition



With an abuse of the language we will say that the above
representation of ~$E$ is a \emph{decomposable} or, respectively,
\emph{indecomposable Horikawa representation}.

\proof
Any (anti-)automorphism of~$Y$ is given by a linear expression in
bihomogeneous coordinates. It is easy to see that an indecomposable
action (whose anti-holomorphic involutions transpose the rulings) can
be converted to a canonical form with $s$ as above and
$\cy\:[\coord(u_1:v_1),\coord(u_2:v_2)]\mapsto
[\coord(\bar u_2:\bar v_2),\coord(\bar u_1:\bar v_1)]$.  A
decomposable action splits into product; up to isomorphism and
interchanging~$\cy$ and $s\circ\cy$ there are two actions on~$\Cp1$:
the holomorphic involution is the map
$s\:\coord(u:v)\mapsto\coord(-u:v)$, and one of the anti-holomorphic
ones is either $c_a\:\coord(u:v)\mapsto\coord(\bar u:\bar v)$ or
$c_b\:\coord(u:v)\mapsto\coord(\bar v:\bar u)$. (Note that
$s\circ c_a$ is isomorphic to~$c_a$, and
$s\circ c_b\:\coord(u:v)\mapsto\coord(-\bar u:\bar v)$ is fixed point
free.) Combining~$c_a$, $c_b$, and $s\circ c_b$, up to permutation of
the factors and of~$\cy$ and~$s\circ\cy$ one obtains the four actions
listed in the statement:
\hbox{\def\lastref{actions}\vbox{\leavevmode\hsize.5\hsize
\roster
\item"\ditto1" $\cy=c_a\times c_a$;
\item"\ditto2" $\cy=c_a\times c_b$;
\endroster}\vbox{\leavevmode\hsize.5\hsize
\roster
\item"\ditto3" $\cy=c_b\times c_b$;
\item"\ditto4" $\cy=c_b\times(s\circ c_b)$.\qed
\endroster}}
\endproof\nofrills

\remark{Remark}
Note that for \iref{actions}1 and~\ditto5
both the expositions have $\er\ne\emptyset$.  Thus, in this paper we
are only concerned with~\iref{actions}2--\ditto4.
\endremark

Homologically, decomposable and indecomposable actions differ by the
induced action in $H_2(Y;\Z)\cong U$: the holomorphic involution
always acts identically, and the anti-holomorphic ones induce
multiplication by~$(-1)$ in the decomposable case and
$y_1\mapsto-y_2$, $y_2\mapsto-y_1$ in the indecomposable one (where
$y_1$, $y_2$ are the classes of the rulings on~$Y$).

\subsection{Special Horikawa models}\label{special}
In the construction of~\ref{3.1} one can as well take for~$Y$ a
rational ruled surface~$\Sigma_2$ with a $(-2)$-section~$S_0$ and
for~$C$, a curve linearly equivalent to a $4$-fold generic section
and disjoint from~$S_0$. ($Y$ can be thought of as the minimal
resolution of the singular point of a quadric cone in~$\Cp3$; then $C$
is cut on~$Y$ by a quartic surface not through the vertex.) Up to
isomorphism there is a unique holomorphic
involution $s\:Y\to Y$ with isolated fixed points.  It
preserves~$S_0$ and has four fixed points:  two in~$S_0$ and two
others in the generatrices through them. As in~\ref{3.1}, if $C$ is
$s$-invariant and does not contain the fixed points of~$s$, precisely
one of the two lifts of~$s$ to~$X$, denoted by~$\tau$, is fixed point
free, and the quotient $E=X/\tau$ is an Enriques surface.  The ruling
of~$Y$ lifts to a $\tau$-invariant pencil in~$X$; its projection
to~$E$ is called the \emph{basic} pencil; its multiple fibers
correspond to the generatrices of~$Y$ through the fixed points
of~$s$. The exceptional section~$S_0$ lifts to two disjoint nodal
curves in~$X$; they are interchanged by~$\tau$ and, thus, project to
one nodal curve in~$E$, called the \emph{basic} nodal curve.

As before, a real structure $\cy\:Y\to Y$ in respect to which
both~$s$ and~$C$ are real defines two real structures on~$E$, called
\emph{expositions} of~$E$. In the $(\Z/2\times\Z/2)$-actions
generated by~$s$ and~$\cy$ one element acts holomorphically and with
isolated fixed points and the two other nontrivial elements act
anti-holomorphically.

\proposition\label{special.actions}
Up to isomorphism there are two such actions\rom:
\roster
\item\local1
with
$\yrr1$ and $\yrr2$ homeomorphic to $S^1\times S^1$ and two invariant
generatrices\rom;
\item\local2
with $\yrr1=S^1\times S^1$, \ $\yrr2=\varnothing$, and no
invariant generatrices.
\endroster
\endproposition

\proof
We blow up the two fixed points of~$s$ not on~$S_0$ and blow down the
generatrices through them. This is an equivariant transformation
whose result is $\Cp1\times\Cp1$ with an action which has invariant
generatrices (say, the image of~$S_0$). The statement follows now
from~\ref{actions}.
\endproof

\subsection{Basic elliptic pencils}\label{3.3}
Consider a (special or nonspecial) real Horikawa representation of a
real Enriques surface $E$.  Denote by $|2P_1|$, $|2P_2|$ the basic
elliptic pencils on~$E$ in the nonspecial case, and by $|2P|$
and~$R$, the basic elliptic pencil and basic nodal curve in the
special case. Let $(P_1,P_1')$, $(P_2,P_2')$, and $(P,P')$ be the
multiple fibers of the pencils. The classes $([P_1],[P_2])\subset L$
and $([P],[P+R])\subset L$ are called the \emph{basic $U$-pairs}; the
sublattice in~$L$ generated by the basic pair is called the
\emph{basic lattice}.  The following is a direct consequence of the
construction:

\lemma\label{U-types}
The basic lattice
of a Horikawa representation
is a $\c_*$-invariant hyperbolic plane
$U\subset L$. If the
representation
is special or decomposable nonspecial,
$\c_*$ acts on~$U$ by multiplication by~$(-1)$. Otherwise
\rom(indecomposable nonspecial representation\rom),
$\c_*$ transposes the
standard generators $[P_1]$ and~$[P_2]$ of~$U$.\qed
\endlemma

\remark{Remark}
Note that basic lattice determines basic $U$-pair. Indeed, $U$ has a
unique pair of standard generators represented by effective divisors.
In the special case one has to distinguish between~$[P]$ and $[P+R]$.
Since $(P+R)\cdot R=-1<0$, this class cannot be represented by a
fiber of an elliptic pencil: such a fiber would contain~$R$ as a
component and, hence, intersect it trivially.
\endremark

A $\c_*$-invariant hyperbolic plane $U\subset L$ is said to be of
\emph{type~$\I$}, or \emph{decomposable}, if $\c_*$ acts on~$U$ as
multiplication by~$(-1)$, and of \emph{type~$\II$}, or
\emph{indecomposable}, if $\c_*$ transposes a pair of generators
of~$U$. In the case of type~$\I$, i.e., $U\subset L^-$, the unordered
pair of values $\delta(x_1)$, $\delta(x_2)$ does not depend on a
choice of a pair $x_1$, $x_2$ of standard generators of~$U$.
According to these values we will further subdivide
type~$\I$ into $\I(0,0)$, $\I(0,w_2)$, and $\I(w_2,w_2)$. The
following is a consequence of Propositions~\ref{actions}
and~\ref{special.actions} (and the fact that always $\delta([R])=0$):

\widestnumber\item{}

\lemma\label{types}
The basic lattice is of type~$\II$ if and only if the representation
is nonspecial indecomposable. In the other cases the types are\rom:
\roster
\item""
$\I(0,0)$
for nonspecial action~\iref{actions}1 and special
action~\iref{special.actions}1\rom;
\item""
$\I(0,w_2)$
for nonspecial action~\iref{actions}2\rom;
\item""
$\I(w_2,w_2)$ for nonspecial actions~\iref{actions}3, \ditto4
and special action~\iref{special.actions}2.\qed
\endroster
\endlemma

\subsection{Selected models of empty surfaces}\label{3.4}
Our
proof  of the main theorem appeals to nonspecial Horikawa
representations
with
decomposable action~\iref{actions}3. As follows from
Lemma~\ref{types} and Remark
in~\ref{3.1}, this type is distinguished by the type of
its
basic lattice, which must be $\I(0,w_2)$. In some affine coordinates
$(x,y)$ on~$Y$ the action is given by $s\:(x,y)\mapsto(-x,-y)$,
$c\:(x,y)\mapsto(\bar x^{-1},\bar y)$. Since $\yrr1=S^1\times S^1$
and $\yrr2=\emptyset$, one has:

\proposition\label{empty}
If the real part of~$E$ is empty, so is~$\crr1$. Conversely, if
$\crr1=\varnothing$, the real part of~$E$ is empty for one and only
one of the expositions.
\endproposition

Let us represent branch curves $C$ by polynomials in $x,y$. Consider
the vector space of $s$-invariant polynomials of bidegree $(4,4)$,
its real part
$$
\tsize
\Cal C=\bigl\{p=\sum a_{i,j}x^iy^j\bigm|i=j\bmod2,\
0\le i,j\le 4,\ \bar a_{4-i,j}=a_{i,j}\bigr\},
$$
the corresponding projective space $P\Cal C$,
and subsets $\Cal M_0\subset\Cal M\subset P\Cal C$
$$
\gather
\Cal M=\bigl\{p\in P\Cal C\bigm|\crr1=\emptyset\bigr\},
\\\noalign{\vskip1\jot}
\Cal M_0=\Bigl\{p\in\Cal M\Bigm|
\vcenter{\offinterlineskip\hbox{%
the curve $C=\{p=0\}$ has at worst simple singularities,}
\kern2pt
\hbox{%
$a_{0,0}\ne 0$, $a_{0,4}\ne0$, $a_{4,0}\ne0$, $a_{4,4}\ne0$}}
\Bigr\}.
\endgather
$$

\theorem\label{model}
For each $p\in\Cal M_0$ the Horikawa construction gives a unique
empty real Enriques surfaces.  The space $\Cal M_0$ is connected\rom;
in particular, all empty real Enriques surfaces obtained by the
Horikawa construction from points of $\Cal M_0$ are of the same
deformation type.
\endtheorem

\proof
The first statement is contained in~\ref{empty}.
The connectedness follows from two observations. First, $\Cal M$ is
convex in $P\Cal C$ and, hence, connected. (Indeed, $\crr1=\emptyset$
if and only if, up to multiplication by~$(-1)$, $p>0$ on~$\yrr1$.)
Second, the conditions $a_{0,0}\ne0$, $a_{0,4}\ne0$, $a_{4,0}\ne0$,
$a_{4,4}\ne0$, as well as the restriction on the singularities, are
of codimension at least~2.  The deformation type of the resulting
surface is preserved since the (unique)
exposition with empty real part,
corresponding to $p>0$ on~$\yrr1$, varies continuously with $p$.
\endproof

\section{Existence of pencils and models}

\subsection{Real elliptic pencils on~$E$}
Recall that an effective divisor $D=\sum m_iR_i$ on~$E$ with
all~$R_i$ distinct and irreducible is called a \emph{divisor of
canonical type} if $D\cdot R_i=K_E\cdot R_i=0$ for all~$i$. It is
called \emph{indecomposable} if it is connected and $\gcd(m_i)=1$.

\lemma\label{D2pencil}
If $D$ is a divisor of canonical type and the class of~$D$ is
primitive in~$L$, then $|2D|$ is an elliptic pencil on~$E$.
\endlemma

\proof
As shown in~\cite{DC}, Proposition~3.1.2, $|D|$ or $|2D|$ is an
elliptic pencil whenever $D$ is an indecomposable divisor of
canonical type.  Let us prove that $D$ is indecomposable. Assume that
$D=\sum n_iD_i$, $n_i\ge1$, with all $D_i$ indecomposable and
$D_i\cdot D_j=0$.  Then $|D_1|$ or $|2D_1|$ is an elliptic pencil,
and the other components are its fibers.  Hence, the classes of
all~$D_i$ are some multiples of the class of a multiple fiber of the
pencil. Since $([D]\bmod\Tors)$ is primitive, this implies that $D=D_1$
and that $|D|$ is not an elliptic pencil.
\endproof

\theorem\label{pencil}
A real Enriques surface $E$ admits a real elliptic pencil if and only
if there exists an element $x\in L^-$ with $x^2=0$. If $E$ is unnodal,
for any primitive element~$x\in L^-$ with $x^2=0$ either~$x$ or~$-x$
can be realized as the class of a multiple fiber of a real elliptic
pencil.
\endtheorem

\proof
The `only if' part is obvious (cf.~\ref{reality}). For the `if' part,
pick a primitive element $x\in L^-$ with $x^2=0$ and an effective
divisor~$D$ with $[D]=\pm x$. According to \cite{DC}, Theorem~3.2.1,
$D\sim D'+\sum m_iR_i$, $m_i\ge1$, where $R_i$ are nodal curves and
$D'$ is a divisor of canonical type.  Moreover, the class $[D']\in L$
is uniquely determined by~$x$ and is obtained from~$x$ by a series of
reflections. Since $x$ is primitive, so is~$[D']$ and, by
Lemma~\ref{D2pencil}, $|2D'|$ is an elliptic pencil. Due to the
uniqueness of~$[D']$ it is $\c$-invariant.
\endproof

\subsection{Existence of real Horikawa models}
The following two pure complex results are known. Their proofs are
found, e.g., in~\cite{BPV}, Theorems~18.1 and~18.2.

\lemma\label{P2map}
If an Enriques surface $E$ has a pair $|2P_1|$, $|2P_2|$ of elliptic
pencils with $P_1\cdot P_2=1$, then $E$ admits a nonspecial Horikawa
representation with $|2P_1|$, $|2P_2|$ as basic pencils.\qed
\endlemma

\lemma\label{PR2map}
If an Enriques surface $E$ has an elliptic pencil $|2P|$ and a nodal
curve ~$R$ with $P\cdot P=1$, then $E$ admits a special Horikawa
representation with $|2P|$ and $R$ as the basic pencil and nodal
curve.\qed
\endlemma

\theorem\label{existence}
A real Enriques surface~$E$ admits a real Horikawa representation if
and only if $L$ contains a $\c_*$-invariant hyperbolic plane~$U$ of
type~$\I$ or~$\II$. If $E$ is unnodal, $U$ can be taken for the basic
lattice of a representation, whose type is determined by
Lemma~\ref{types}.  In general, if $U$ is of type~$\I$
\rom(respectively, type~$\II$\rom), $E$ admits either a special
Horikawa representation or a nonspecial decomposable
\rom(respectively, indecomposable\rom)  Horikawa representation.
\endtheorem

\proof
The necessity of the condition follows from~\ref{U-types}. To prove
the `if' part, consider the standard generators $x_1,x_2\in U$
realized by effective divisors. According to \cite{DC}, Lemma~3.3.1,
by a sequence of reflections against classes of nodal curves $x_1$,
$x_2$ can be taken to some classes $y_1,y_2\in L$, unique up to
reordering, so that either
\roster
\item\local1
$y_1=[D_1]$ and $y_2=[D_2]$, or
\item\local2
$y_1=[D_1]$ and $y_2=[D_1+R]$,
\endroster
where $D_1$, $D_2$ are indecomposable divisors of canonical type and
$R$ is a nodal curve. In case~\loccit1 the order of~$y_1$, $y_2$ is
also determined by~$x_1$, $x_2$; hence, the pair $(y_1,y_2)$ behaves
in respect to~$\c_*$ in the same way as $(x_1,x_2)$ and a real
Horikawa representation is given by Lemma~\ref{P2map}. In
case~\loccit2 the uniqueness implies that $y_1$, $y_2$ are
$\c_*$-skew-invariant and a representation is given by
Lemma~\ref{PR2map}.  Finally, if $E$ is unnodal, no reflection is
necessary and $x_1$, $x_2$ are themselves realized by indecomposable
divisors of canonical type.
\endproof

\section{Proof of the Main Theorem}\label{s6}

\subsection{Eigenlattices of~$\c_*$}
Let $D_4$ denote the negative lattice generated by the root system of
the same name.

\lemma\label{lattices}
For an empty real Enriques surface $L^+ = D_4$ and
$L^-= D_4\oplus U$.
\endlemma

\proof
It suffices to prove that $L^+=D_4$, since $D_4$ admits a unique,
up to isomorphism, primitive embedding into $E_8\oplus U$,
see~\cite{Nik79}, Theorem 1.14.2 and Remark 1.14.6.

Since $E_\R$ is empty, $\rank L^+=4$, $\sigma(L^+)=-4$, and
the discriminant form $\discr L^+$ of~$L^+$
is even.  The dimension of $\discr L^+$ is either~$0$, or~$2$,
or~$4$.  In the first case $L^+$ is unimodular and must have
signature $0\bmod8$ (see~\cite{Serre}). In the last case
$L^+(\frac12)$ is unimodular and even, which also contradicts to the
signature congruence.  Thus, $L^+$ is an even integral negative
lattice of rank~$4$ with $\discr L^+$ even and of
dimension~$2$.  The only such lattice is $D_4$, see~\cite{Nik79},
Remark 1.14.6.
\endproof

\remark{Remark}
Our case is very special and much simpler
than the general situation treated by Nikulin. For example, the fact that
$D_4$ is determined up to isomorphism by its genus is
proved in a few lines. Indeed, since the
discriminant group contains a vector of square~$1$, any such lattice
is a sublattice of index~$2$ of a negative unimodular lattice of
rank~$4$, i.e., $4\langle-1\rangle$. Since the latter is odd, the
lattice in question is uniquely determined as its maximal even
sublattice. For a similar reason the orthogonal complement of~$D_4$
in $E_8\oplus U$ is also unique, and the uniqueness of an embedding
$D_4\subset E_8\oplus U$ is then straightforward.
\endremark


\lemma\label{pair}
If $E$ is an empty real Enriques surface, $L$ contains a hyperbolic
plane of type $\I(0,w_2)$.
\endlemma

\proof
Due to~\ref{lattices} $L^-\cong D_4\oplus U$ does contain a
hyperbolic plane~$U$. Since $U$ is unimodular and $D_4$ is determined
by its genus, the embedding $U\to L^-$ is unique up to isomorphism.
Let $x_1$, $x_2$ be a standard pair of generators of~$U$, and
$e_1,\dots,e_4$ some standard generators of~$D_4$ (so that
$e_i^2=-2$, $e_ie_j=0$ for $1\le i<j\le3$, and $e_ie_4=1$ for
$1\le i\le3$). Consider the following three possibilities:

\def\Case#1{\par\leavevmode\hbox{{\it Case #1\/}:\enspace}%
\ignorespaces}

\smallskip
\Case1
$\delta(x_1)=0$ and $\delta(x_2)=w_2$. Then $U$ is of
type~$\I(0,w_2)$.

\Case2
$\delta(x_1)=\delta(x_2)=w_2$. If there is an index~$i$ with
$\delta(e_i)=0$, the sublattice generated by $x_1+x_2+e_i$ and $x_1$
is of type $\I(0,w_2)$. Otherwise, $\delta(e_1+e_4)=0$ and the
sublattice generated by $x_1+x_2+e_1+e_4$ and~$x_1$ is of type
$\I(0,w_2)$.

\Case3
$\delta(x_1)=\delta(x_2)=0$ for {\it any\/} hyperbolic plane
$U\subset L^-$. Take for~$U$ the basic lattice of a real Horikawa
representation, which exists due to Theorem~\ref{existence}. If the
representation is nonspecial and decomposable, its basic pencils
$|2P_1|$, $|2P_2|$ have real multiple fibers and, hence, the only
intersection point of~$P_1$ and~$P_2$ is real. If the representation
is special, the basic pencil~$|2P|$ has real multiple fibers and
$x_2-x_1$ is realized by a real nodal curve~$R$, and the intersection
point of~$P$ and~$R$ is real. Existence of a real point contradicts
to the assumption $E_\R=\varnothing$.
\endproof

\subsection{Proof of the theorem}\label{last}

\lemma\label{gen}
Any real Enriques surface can be made unnodal by a small real
deformation.
\endlemma

\proof
Surfaces with nodal curves form a countable union of hypersurfaces in
the period space of complex Enriques surfaces: these hypersurfaces
are the hyperplane sections defined by $(-2)$-vectors in the
$K3$-lattice (see, e.g., the construction of the period space
in~\cite{BPV}). For our purpose, it is sufficient to consider the
period space of local deformations. Due to the local Torelli theorem,
the real part of the local period space is nonsingular and represents
small real deformations of a given real Enriques surface
(cf.~\cite{Kh76}).  Thus, in any neighborhood of a real point there
are real points representing unnodal surfaces.
\endproof

\lemma\label{M}
Any unnodal empty real Enriques surface is obtained by the
Ho\-ri\-ka\-wa construction from a point of~$\Cal M_0$.
\endlemma

\proof
Due to~\ref{pair} $L$ contains a hyperbolic plane of
type~$\I(0,w_2)$, which, due to~\ref{existence} (and since $E$ is
unnodal) is the basic lattice of a real Horikawa representation based
on a decomposable action.  Due to~\ref{types}, this is the action
selected in~\ref{3.4}.
\endproof

\proof[Proof of Theorem~\ref{main}]
To deform one empty real Enriques surface to another we deform them
both to unnodal surfaces (see~\ref{gen}), which, due to~\ref{M}, are
represented by points of $\Cal M_0$, and connect the results by a
path in~$\Cal M_0$ (see~\ref{model}).
\endproof

\remark{Remark}
The choice of one of the two lifts of the real structure on an empty
real Enriques surface to the covering $K3$-surface
defines a nonramified double covering of the moduli space of empty
real Enriques surfaces.
This covering is nontrivial:
the two covering real structures can be exchanged by a
diffeomorphism or by the monodromy along a loop in the moduli
space. This is easily seen on Horikawa models
with decomposable action ~\iref{actions}2 or ~\iref{actions}4. In
particular,
similar to empty real Enriques surfaces, the $K3$-surfaces equipped
with an ordered pair of commuting anti-holomorphic fixed point free
involutions
with fixed point free composition are all of the same deformation type.
\endremark

%
%
\enddocument